\documentclass{amsart}
\usepackage{diagrams}
\theoremstyle{plain}
 \newtheorem{thm}{Theorem}[section]
 \numberwithin{equation}{section} %% Comment out for sequentially-numbered
 \numberwithin{figure}{section} %% Comment out for sequentially-numbered
 \theoremstyle{plain}
 \newtheorem{prop}[thm]{Proposition} %%Delete [thm] to re-start numbering
\theoremstyle{plain}
\newtheorem{lemma}[thm]{Lemma}
 \theoremstyle{plain}
 \newtheorem{cor}[thm]{Corollary}
 \theoremstyle{definition}
 \newtheorem{defn}[thm]{Definition}
 \theoremstyle{definition}
  \newtheorem{example}[thm]{Example}
 \theoremstyle{remark}
 \newtheorem{rem}[thm]{Remark}
 \theoremstyle{remark}
 \newtheorem*{note*}{Note}
 \theoremstyle{remark}
 
 \theoremstyle{remark}
 \newtheorem*{acknowledgement*}{Acknowledgement}

\title{On Weak Lie 2-algebras}
\author{Dmitry Roytenberg}
\address{IH\'{E}S, Bures-sur-Yvette, France}
\email{roytenbe@ihes.fr} \keywords{Categorification, Lie algebra,
crossed module, Courant algebroid}

\subjclass{Primary 17B55; Secondary 18G55, 58H15, 55U15, 17B81}

\begin{document}

\begin{abstract} A Lie 2-algebra
is a linear category equipped with a functorial bilinear operation
satisfying skew-symmetry and Jacobi identity up to natural
transformations which themselves obey coherence laws of their own.
Functors and natural transformations between Lie 2-algebras can also
be defined, yielding a 2-category. Passing to the normalized chain
complex gives an equivalence of 2-categories between Lie 2-algebras
and 2-term "homotopy everything" Lie algebras; for strictly
skew-symmetric Lie 2-algebras, these reduce to $L_\infty$-algebras,
by a result of Baez and Crans. Lie 2-algebras appear naturally as
infinitesimal symmetries of solutions of the Maurer--Cartan equation
in some differential graded Lie algebras and $L_\infty$-algebras. In
particular, (quasi-) Poisson manifolds, (quasi-) Lie bialgebroids
and Courant algebroids provide large classes of examples.
\end{abstract}

\maketitle
\section{Introduction}
The purpose of this note is to complete the categorification of the
notion of Lie algebra started by Baez and Crans in \cite{BaezCrans},
who introduced \emph{semi-strict Lie 2-algebras}. By definition,
such a structure is given by a bilinear bracket operation on a
linear category, which is strictly skew-symmetric but obeys the
Jacobi identity only up to a coherent trilinear natural
transformation, called the Jacobiator. It was further shown that,
upon passing to the normalized chain complex, such a structure is
equivalent to a 2-term $L_\infty$-algebra. Passing to cohomology and
using homotopy invariance, one gets a classification of semistrict
Lie 2-algebras in terms of 3rd Chevalley-Eilenberg cohomology of a
Lie algebra with coefficients in a module.

From the point of view of category theory, this picture is somewhat
incomplete, as the skew-symmetry holds strictly as an equation.
Besides, there exist structures -- Leibniz algebras and Courant
algebroids -- where the Jacobi holds strictly (in the form of a
Leibniz identity), whereas it is the skew-symmetry that is weakened.
It was shown in \cite{RoyWe} for the case of Courant algebroids that
skew-symmetrizing the bracket does lead to a 2-term
$L_\infty$-algebra; however, some information is lost in the
process, and besides, it is not clear why this has to happen,
conceptually. That is why, to have a better understanding of this
phenomenon, we propose to work from the outset with \emph{weak Lie
2-algebras}, where both the skew-symmetry and the Jacobi identity
are allowed to hold only up to coherent natural transformations --
the \emph{alternator} and the \emph{Jacobiator}. This structure
resembles (and perhaps is Koszul dual to, in a sense yet to be made
preceise) a linear braided monoidal category. Passing to the
normalized chain complex leads to a new structure -- \emph{
$EL_\infty$-algebra} ($E$ for (homotopy) "everything"). Weak Lie
2-algebras form a 2-category, and so do 2-term $EL_\infty$-algebras,
and these 2-categories are equivalent (Theorem \ref{thm:2catequiv}).
We claim that weak Lie 2-algebras are the correct categorification
of Lie algebras and henceforth refer to them simply as Lie
2-algebras. Lie 2-algebras with trivial alternator will be referred
to as semistrict, as in \cite{BaezCrans}; those with trivial
Jacobiator -- hemistrict. These form full sub-2-categories. This is
the content of the first section.

In the next section, we describe the skew-symmetrization functor,
which is a projection onto the sub-2-category of semistrict Lie
2-algebras (Theorem \ref{thm:skewsym}). There is some "fudging"
involved here: the usual inverse factorials must be multiplied in
some places by certain additional rational numbers in order for the
Theorem to hold. The origin and meaning of these factors is, at this
point, unclear.

In the following section, we discuss homotopy invariance of Lie
2-algebras and deduce a classification of skeletal Lie 2-algebras
(Theorems \ref{thm:invariance} and \ref{thm:classification}). The
classification uses what appears to be a new cohomology theory for
Lie algebras, based on the fact that a Lie algebra is a Leibniz
algebra which is also skew-symmetric: the Jacobiator defines a
Loday-Pirashvili 3-cocycle whose behavior under permutations of the
arguments is controlled by the alternator, and the fact that
coboundaries are cocycles depends both on the skew-symmetry and
Jacobi identity of the Lie algebra. This is similar to Eilenberg and
MacLane's cohomology theory for abelian groups which uses both the
associativity and commutativity of the group law.
Skew-symmetrization induces a map from the new cohomology onto the
3rd Chevalley-Eilenberg cohomology (Theorem \ref{thm:LPtoChE}); this
map is an isomorphism if, and only if, the alternator is symmetric.

Finally, in the last section, we discuss applications of the theory
to questions in deformation theory. Namely, given a solution of the
Maurer-Cartan equation in a differential graded Lie algebra
concentrated in degrees $(-3,+\infty)$, we construct the hemistrict
Lie 2-algebra of its inner symmetries, mapping to the (ordinary) Lie
algebra of infinitesimal symmetries, forming a categorified crossed
module (Theorem \ref{thm:2crossed}). Such a dgla controls, for
instance, the deformation theory of Courant algebroids
\cite{Roy4-GrSymp} and hence, this theorem, combined with Theorem
\ref{thm:skewsym}, generalizes the main result of \cite{RoyWe}. The
construction itself uses \emph{derived brackets}, as described in
\cite{KS5}; in case of dgla's concentrated in degrees
$(-2,+\infty)$, a similar but simpler construction is well known and
yields the crossed module of infinitesimal symmetries.

\begin{acknowledgement*} The bulk of the research for this paper was
carried out at the Max Planck Institut f{\"{u}}r Mathematik in Bonn
where I spent the months of February through June of 2007. I would
like to thank the Institute for providing excellent working
conditions as well as financial support. My gratitude goes out also
to the organizers of the XXVI Workshop on Geometric Methods in
Physics in Bialowieza (Poland), and of the Program on Poisson Sigma
models, Lie algebroids, deformations and higher analogues at the
Erwin Schr{\"{o}}dinger Institute in Vienna, for giving me the
opportunity to present the results of this paper. The text of the
paper was written up at IH\'{E}S.
\end{acknowledgement*}

\section{Lie 2-algebras and $EL_\infty$-algebras}\label{sec:Lie2}

\subsection{Categorified linear algebra}
Fix a ground field $\mathbf{k}$. It can be assumed arbitrary, except
where indicated. We let $\mathbf{Vect}$ denote the category of
vector spaces over $\mathbf{k}$. In what follows we shall freely use
the notations and terminology of \cite{BaezCrans}, with minor
modifications. Thus, a \emph{2-vector space} is a linear category,
i.e. a category internal to $\mathbf{Vect}$: objects and morphisms
form vector spaces, and all the structure maps are linear. If $V$ is
such a category, we denote its space of objects by $V_0$, its space
of arrows by $V_1$, the source and target maps, respectively, by
$s,t:V_1\rightarrow V_0$, the identity map by $1:V_0\rightarrow V_1$
(with $x\mapsto 1_x$), and the composite of the arrows
$\mathbf{a}:x\rightarrow y$ and $\mathbf{b}:y\rightarrow z$ by
$\mathbf{ba}:x\rightarrow z$. It is clear what it means for a
functor between two linear categories, or a natural transformation
between two such functors, to be linear. Given two linear functors
$F:U\rightarrow V$ and $G:V\rightarrow W$, we denote their
composition by $G*F:U\rightarrow W$. Given two linear natural
transformations $\Phi:F\Rightarrow F'$ and $\Psi:G\Rightarrow
G'$, with $F,F':U\rightarrow V$ and $G,G':V\rightarrow W$, we denote
their horizontal composite by $\Psi*\Phi:G*F\Rightarrow G'*F'$;
whiskering of a natural transformation by a functor is also denoted
by $*$. Finally, given natural transformations $\Phi:F\Rightarrow
G$ and $\Psi:G\Rightarrow H$, we denote their vertical composite by
$\Psi\circ\Phi:F\Rightarrow H$. The definitions of these
compositions are standard in category theory, and it is trivial to
check that they preserve linearity. Thus we have the (strict)
2-category $\mathbf{2Vect}$.

The standard constructions of linear algebra carry over to
$\mathbf{2Vect}$ in an obvious manner. In particular, the tensor
product of 2-vector spaces is defined "dimension-wise" and satisfies
the usual universal property with respect to multilinear functors.
This makes $\mathbf{2Vect}$ a symmetric monoidal 2-category; we
denote the action of the transposition by $\sigma:V\otimes
W\rightarrow W\otimes V$. In fact, $\mathbf{2Vect}$ is a
\emph{closed} symmetric monoidal 2-category which means, in
particular, that $\mathrm{Hom}_{\mathbf{2Vect}}(V,W)$ is naturally a
linear category, and the composition of linear functors is bilinear.

Given a 2-vector space $V$, we define its \emph{normalized cochain
complex} $N(V)$ by
\begin{eqnarray*}
    N(V)^0 &=& V_0 \\
  N(V)^{-1}&=& \ker(s)
\end{eqnarray*}
with $d:N(V)^{-1}\rightarrow N(V)^0$ given by the restriction of
$t$. It is easy to check (see \cite{BaezCrans}) that this extends to
a (strict) 2-functor from $\mathbf{2Vect}$ to the 2-category
$\mathbf{2Term}$, consisting of 2-term cochain complexes
(concentrated in degrees $(-2,0]$), chain maps and chain homotopies.
In fact, the normalization functor has a quasi-inverse $\Gamma$,
given on objects by
\begin{eqnarray*}
  \Gamma(C)_0 &=& C^0 \\
  \Gamma(C)_1 &=& C^0\oplus C^{-1}
\end{eqnarray*}
with
\begin{eqnarray*}
  s(x,a) &=& x\\
  t(x,a) &=& x+da \\
  1_x &=& (x,0) \\
  (y,b)(x,a) &=& (x,a+b)\qquad \mathrm{if}\quad y=x+da
\end{eqnarray*}

\begin{thm} (\cite{BaezCrans})
The 2-functors $N$ and $\Gamma$ give
an equivalence of 2-categories $\mathbf{2Vect}$ and
$\mathbf{2Term}$.
\end{thm}

\begin{rem}
The functors $N$ and $\Gamma$ can be defined in a more general
setting of simplicial vector spaces and non-positively graded
cochain complexes. That they define an equivalence of these two
categories is a classical theorem in homological algebra, due to
Dold and Kan. The categorical equivalence in the theorem above
follows from this by applying the nerve functor from linear
categories to simplicial vector spaces. The souped-up 2-categorical
version is due to Baez and Crans.
\end{rem}

\begin{rem}
The simple observation underlying the above result is that any arrow
$\mathbf{a}:x\rightarrow y$ in a linear category can be uniquely
decomposed as $\mathbf{a}=1_x+a$, where $a=\mathbf{a}-1_x\in\ker{s}$
is called the \emph{arrow part} of $\mathbf{a}$. The linearity of
composition then forces it to be just the addition of arrow parts.
It also implies that any linear category is in fact a
\emph{groupoid}, for $(x+da,-a)(x,a)=(x,0)$.
\end{rem}

In the sequel, we shall freely use the canonical isomorphism
$V\simeq\Gamma N(V)$ and write $\mathbf{a}=x+a$ when it is not
likely to cause confusion.

\subsection{Multilinear operations} Using the Dold-Kan
correspondence, it is tempting to conclude right away that
multilinear functorial operations on 2-vector spaces are in
one-to-one correspondence with multilinear chain maps on their
normalized complexes. The problem is that the normalization functor
does not commute with tensor products; in fact, $\mathbf{2Term}$ is
not even closed under the tensor product of cochain complexes. A
proper treatment of this problem (certainly necessary if we want
eventually to understand higher linear categories) requires a
careful analysis of the behavior of the nerve and normalization
functors with respect to tensor products and the Eilenberg-Zilber
construction. This issue will be addressed elsewhere. Here we
provide instead a quick fix based on the following observation.

\begin{prop}\label{prop:derived}
Let $V_1,\ldots,V_n,V$ be linear categories and
$T:V_1\otimes\cdots\otimes V_n\rightarrow V$ a linear functor. Then
\begin{enumerate}
  \item $\forall a_k\in N(V_k)^{-1}$ and $x_i\in N(V_i)^0$,
  $i\neq k$, $T(x_1,\ldots,x_{k-1},a_k,x_{k+1},\ldots,x_n)\in
  N(V)^{-1}$ and
  $$dT(x_1,\ldots,x_{k-1},a_k,x_{k+1},\ldots,x_n)=T(x_1,\ldots,x_{k-1},da_k,x_{k+1},\ldots,x_n)$$
  \item For $a_i\in N(V_i)^{-1}$, $a_j\in N(V_j)^{-1}$ and
  arbitrary other arguments,
  $$T(\ldots,a_i,\ldots,a_j,\ldots)=T(\ldots,da_i,\ldots,a_j,\ldots)=T(\ldots,a_i,\ldots,da_j,\ldots)$$
\end{enumerate}
\end{prop}

It follows that $T$ is completely determined by its value on objects
and on arrows of the form $x_1\otimes\cdots\otimes x_{k-1}\otimes
a_k\otimes x_{k+1}\otimes\cdots\otimes x_n$. More precisely,

\begin{cor}\label{cor:normequiv}
The linear categories
$\mathrm{Hom}_{\mathbf{2Vect}}(V_1\otimes\cdots\otimes V_n,V)$ and\\
$\mathrm{Hom}_{\mathbf{Ch}}(N(V_1)\otimes\cdots\otimes
N(V_n),N(V))$, where $\mathbf{Ch}$ is the 2-category of
non-positively graded cochain complexes, are canonically isomorphic.
\end{cor}

Let us spell out what the above proposition says in the case of a
binary operation on a linear category $V$:

\begin{prop}\label{prop:crossedmodule}
Let $[\cdot,\cdot]:V\otimes V\rightarrow V$ be a bilinear functor.
Then
\begin{equation}\label{eqn:catbracket}
    [(x,a),(y,b)]=([x,y],[x,b]+[a,y]+[a,b])
\end{equation}
and the following \emph{crossed module identities} hold:
\begin{eqnarray}
  d[x,b] &=&[x,db]\label{eqn:dg1} \\
  d[a,y] &=&[da,y]\label{eqn:dg2} \\
  {[}da,b] &=&[a,b]=[a,db]\\
  d[a,b] &=&[da,db]
\end{eqnarray}
\end{prop}

The corresponding bracket $[\cdot,\cdot]:N(V)\otimes N(V)\rightarrow
N(V)$, given by
\begin{equation}\label{eqn:chainbracket}
    [(x,a),(y,b)]=([x,y],[x,b]+[a,y]),
\end{equation}
is then a chain map. Conversely, any such operation on a 2-term
chain complex $C$ uniquely determines a bilinear functorial bracket
on $\Gamma(C)$ by setting $[a,b]$ to be the \emph{derived bracket}:
\begin{equation}\label{eqn:derivedbracket}
    [a,b]=[da,b]=[a,db]
\end{equation}

\subsection{The 2-category of Lie 2-algebras}
We are now ready to define weak Lie 2-algebras.

\begin{defn}
A \emph{Lie 2-algebra} is a linear category $L$ equipped with the
following structure:
\begin{itemize}
\item a bilinear functor $[\cdot,\cdot]:L\otimes L\rightarrow L$,
called the \emph{bracket};
\item a bilinear natural transformation
$$S:[\cdot,\cdot]\Rightarrow-[\cdot,\cdot]*\sigma,$$ called the
\emph{alternator};
\item a trilinear natural transformation
$$J:[\cdot,[\cdot,\cdot]]\Rightarrow[[\cdot,\cdot],\cdot]+[\cdot,[\cdot,\cdot]]*\sigma_{12},$$
called the \emph{Jacobiator}.
\end{itemize}
In addition, the following diagrams are required to commute:

%Jacobiator hexagon
\def\OT{[x,[y,[z,w]]]}
\def\OLT{[x,[[y,z],w]]+[x,[z,[y,w]]]}
\def\ORT{[[x,y],[z,w]]+[y,[x,[z,w]]]}
\def\OLB{\begin{array}{c}
           [[x,[y,z]],w]+[[y,z],[x,w]]+ \\
           +[[x,z],[y,w]]+[z,[x,[y,w]]]
         \end{array}
}
\def\ORB{\begin{array}{c}
           [[[x,y],z],w]+[z,[[x,y],w]]+ \\
           +[y,[[x,z],w]]+[y,[z,[x,w]]]
         \end{array}
}
\def\OB{\begin{array}{c}
          [[[x,y],z],w]+[[y,[x,z]],w]+[[y,z],[x,w]]+\\
          +[[x,z],[y,w]]+[z,[[x,y],w]]+[z,[y,[x,w]]]
        \end{array}
}
\def\ALT{{[x,J_{y,z,w}]}}
\def\ART{J_{x,y,[z,w]}}
\def\ALM{J_{x,[y,z],w}+J_{x,z,[y,w]}}
\def\ARM{J_{[x,y],z,w}+[y,J_{x,z,w}]}
\def\ALB{[J_{x,y,z},w]+1+1+[z,J_{x,y,w}]}
\def\ARB{1+1+J_{y,[x,z],w}+J_{y,z,[x,w]}}

\begin{diagram}[labelstyle=\scriptstyle]
&&\OT&&\\
&\ldTo[nohug](2,2)<{\ALT}&&\rdTo[nohug](2,2)>{\ART}&\\
\OLT&&&&\ORT\\
\dTo<{\ALM}&&&&\dTo>{\ARM}\\
\OLB&&&&\ORB\\
&\rdTo[nohug](2,2)<{\ALB}&&\ldTo[nohug](2,2)>{\ARB}&\\
&&\OB&&\\
\end{diagram}

%alternator triangle
\begin{diagram}
[[x,y],z]&                                   &\rTo^{[S_{x,y},z]} &                                    &-[[y,x],z]\\
         &\rdTo[nohug](2,2)<{\hat{J}_{x,y,z}}&                   &\ldTo[nohug](2,2)>{-\hat{J}_{y,x,z}}&          \\
         &                                   &[x,[y,z]]-[y,[x,z]]&                                    &          \\
\end{diagram}

%alternator square
\begin{diagram}
          [x,[y,z]]&   \rTo^{[x,S_{y,z}]}         &-[x,[z,y]]          \\
   \dTo<{J_{x,y,z}}&                              &\dTo>{-J_{x,z,y}}    \\
[[x,y],z]+[y,[x,z]]&\rTo^{S_{[x,y],z}+S_{y,[x,z]}}&-[[x,z],y]-[z,[x,y]]\\
\end{diagram}

%alternator pentagon
%\begin{diagram}
%                   &                                           &[x,[y,z]]& \rTo^{S_{x,[y,z]}} &-[[y,z],x]&                                    &                    \\
%                   &\ldTo[nohug](2,2)<{J_{x,y,z}}              &         &                    &          &\rdTo[nohug](2,2)>{-\hat{J}_{y,z,x}}&                    \\
%[[x,y],z]+[y,[x,z]]&                                           &         &                    &          &                                    &-[y,[z,x]]+[z,[y,x]]\\
%                   &\rdTo[nohug](3,2)<{[S_{x,y},z]+[y,S_{x,z}]}&         &                    &          & \ldTo[nohug](3,2)>{1+S_{z,[y,x]}}  &                    \\
%                   &                                           &         &-[[y,x],z]-[y,[z,x]]&          &                                    &                    \\
%\end{diagram}

%alternator quasi-symmetry
\begin{diagram}
[x,[y,z]]&                               &\rTo^{1_{[x,[y,z]]}}&                                &[x,[y,z]]\\
         &\rdTo[nohug](2,2)<{S_{x,[y,z]}}&                  &\ldTo[nohug](2,2)>{-S_{[y,z],x}}&         \\
         &                               &-[[y,z],x]        &                                &         \\
\end{diagram}
\end{defn}

\begin{rem}
The natural transformation $\hat{J}$ appearing in the diagrams is,
essentially, the inverse of the Jacobiator:
$$\hat{J}_{x,y,z}=J^{-1}_{x,y,z}-1_{[y,[x,z]]}:[[x,y],z]\rightarrow[x,[y,z]]-[y,[x,z]].$$
It carries the same information as $J$, but with a slight shift in
emphasis: while $J$ measures the failure of $ad(x)=[x,\cdot]$ to be
a derivation of $[\cdot,\cdot]$, $\hat{J}$ measures the failure of
$ad$ to send $[\cdot,\cdot]$ to the commutator bracket of
endomorphisms.
\end{rem}

\begin{rem}
Notice that we do not require $S$ to be symmetric in the sense that
$-S_{y,x}S_{x,y}=1_{[x,y]}$; instead, we impose a weaker condition
$$-S_{[y,z],x}S_{x,[y,z]}=1_{[x,[y,z]]}$$
(the last triangle above). It is natural to wonder whether even this
weakened symmetry assumption can be avoided, but it appears to be
necessary for Theorem \ref{thm:skewsym} to hold. In all examples we
consider $S$ is, in fact, symmetric.
\end{rem}

\begin{defn}
A Lie 2-algebra $L$ is called
\begin{itemize}
\item \emph{semistrict} if $S=1$;
\item \emph{hemistrict} if $J=1$;
\item \emph{strict} if it is both hemistrict and semistrict
\end{itemize}
\end{defn}

\begin{defn}
A morphism of Lie 2-algebras from $(L,[\cdot,\cdot],S,J)$ to
$(L',[\cdot,\cdot]',S',J')$ consists of:
\begin{itemize}
\item  a linear functor $$F:L\rightarrow L'$$
\item  a linear natural transformation
$$F^2:[\cdot,\cdot]'*(F\otimes F)\Rightarrow F*[\cdot,\cdot]$$
such that the following diagrams commute:

%Functor square
\begin{diagram}
[F(x),F(y)]'         &\rTo^{F^2_{x,y}} &F([x,y])         \\
\dTo<{S'_{F(x),F(y)}}&                 &\dTo>{F(S_{x,y})}\\
-[F(y),F(x)]'        &\rTo^{-F^2_{y,x}}&-F([y,x])        \\
\end{diagram}

%Functor hexagon
\begin{diagram}
 [F(x),[F(y),F(z)]']'&\rTo^{J'_{F(x),F(y),F(z)}}&[[F(x),F(y)]',F(z)]'+[F(y),[F(x),F(z)]']'\\
\dTo<{[1,F^2_{y,z}]'}&                          &\dTo>{[F^2_{x,y},1]'+[1,F^2_{x,z}]'}     \\
     [F(x),F([y,z])]'&                          &[F([x,y]),F(z)]'+[F(y),F([x,z])]'        \\
 \dTo<{F^2_{x,[y,z]}}&                          &\dTo>{F^2_{[x,y],z}+F^2_{y,[x,z]}}       \\
         F([x,[y,z]])&   \rTo^{F(J_{x,y,z})}    &F([[x,y],z]+[y,[x,z]])                   \\
\end{diagram}
\end{itemize}
\end{defn}

\begin{defn}
Given two morphisms $(F,F^2):L\rightarrow L'$ and
$(G,G^2):L'\rightarrow L''$, their composite is defined to be
$(G*F,(G*F)^2)$, where
$$(G*F)^2=(G*F^2)\circ(G^2*(F\otimes F))$$
\end{defn}

\begin{rem}
The definition of $(G*F)^2$ is best understood as the total
composite of the following pasting diagram:

%Composition total
\begin{diagram}
          L         &                \rTo^{F}              &          L'         &                \rTo^{G}              &         L''          \\
\uTo<{[\cdot,\cdot]}&\luImplies[nohug][shortfall=1cm]>{F^2}&\uTo~{[\cdot,\cdot]'}&\luImplies[nohug][shortfall=1cm]>{G^2}&\uTo>{[\cdot,\cdot]''}\\
      L\otimes L    &          \rTo_{F\otimes F}           &     L'\otimes L'    &          \rTo_{G\otimes G}           &     L''\otimes L''   \\
\end{diagram}

which can be broken down as the vertical composite of the following
two:

\pagebreak
%Composition Top
\begin{diagram}
          &                                 &                 L              &                                  &  &&\\
          &\ruTo[nohug](2,2)<{[\cdot,\cdot]}&                                &     \rdTo[nohug](2,2)>{F}        &  &&\\
L\otimes L&                                 &\uImplies[shortfall=1.5cm]>{F^2}&                                  &L'&\rTo^{G}&L''\\
          & \rdTo[nohug](2,2)<{F\otimes F}  &                                &\ruTo[nohug](2,2)>{[\cdot,\cdot]'}&  &&\\
          &                                 &           L'\otimes L'         &                                  &  &&\\
\end{diagram}

%Composition Bottom
\begin{diagram}
          &                 &            &                                  &               L'               &                                   &\\
          &                 &            &\ruTo[nohug](2,2)<{[\cdot,\cdot]'}&                                &       \rdTo[nohug](2,2)>{G}       &\\
L\otimes L&\rTo^{F\otimes F}&L'\otimes L'&                                  &\uImplies[shortfall=1.5cm]>{G^2}&                                   &L''\\
          &                 &            &  \rdTo[nohug](2,2)<{G\otimes G}  &                                &\ruTo[nohug](2,2)>{[\cdot,\cdot]''}&\\
          &                 &            &                                  &          L''\otimes L''        &                                   &\\
\end{diagram}

\end{rem}

\begin{defn}
Given two morphisms $(F,F^2),(G,G^2):L\rightarrow L'$, a 2-morphism
$\Theta:(F,F^2)\Rightarrow(G,G^2)$ is a linear natural
transformation $\Theta:F\Rightarrow G$ making the following diagram
commute:
%Natural transformation
\begin{diagram}
               [F(x),F(y)]'&\rTo^{F^2_{x,y}}&F([x,y])             \\
\dTo<{[\Theta_x,\Theta_y]'}&                &\dTo>{\Theta_{[x,y]}}\\
               [G(x),G(y)]'&\rTo^{G^2_{x,y}}&G([x,y])             \\
\end{diagram}
\end{defn}

The horizontal and vertical composites of 2-morphisms are defined to
be those of the corresponding natural transformations. With these
definitions, it is now a matter of routine verification to obtain

\begin{prop}
Lie 2-algebras, their morphisms and 2-morphisms form a strict
2-category, denoted $\mathbf{2Lie}$.
\end{prop}

\subsection{The 2-category of 2-term $EL_\infty$-algebras.}

Let us now apply the normalization functor to the above
construction. Given a Lie 2-algebra $(L,[\cdot,\cdot],S,J)$, denote
the normalized cochain complex $N(L)$ of the underlying linear
category $L$ by $C=C^{-1}\stackrel{d}{\longrightarrow}C^0$. The
induced bracket on $C$ is given by the formula
\ref{eqn:chainbracket} and is a chain map from $C\otimes C$ to $C$.
Furthermore, writing
\begin{eqnarray}
    S_{x,y} &=& ([x,y],-\langle x,y\rangle) \label{eqn:S}\\
  J_{x,y,z} &=& ([x,[y,z]],-\langle x,y,z\rangle) \label{eqn:J}
\end{eqnarray}
and using the naturality of $S$ and $J$, we obtain
\begin{eqnarray}
  {[}x,y]+[y,x] &=& d\langle x,y\rangle\label{eqn:sym1} \\
  {[}a,y]+[y,a] &=& \langle da,y\rangle\label{eqn:sym2}\\
  {[}x,b]+[b,x] &=& \langle x,db\rangle\label{eqn:sym3}
\end{eqnarray}
and
\begin{eqnarray}
  {[}x,[y,z]]-[[x,y],z]-[y,[x,z]] &=& d\langle x,y,z\rangle \label{eqn:Leib1}\\
  {[}a,[y,z]]-[[a,y],z]-[y,[a,z]] &=& \langle da,y,z\rangle \label{eqn:Leib2}\\
  {[}x,[b,z]]-[[x,b],z]-[b,[x,z]] &=& \langle x,db,z\rangle \label{eqn:Leib3}\\
  {[}x,[y,c]]-[[x,y],c]-[y,[x,c]] &=& \langle x,y,dc\rangle \label{eqn:Leib4}
\end{eqnarray}
In other words, $\langle\cdot,\cdot\rangle:C\otimes C\rightarrow
C[-1]$ is a chain homotopy
\begin{equation*}
    \langle\cdot,\cdot\rangle:[\cdot,\cdot]+[\cdot,\cdot]*\sigma\Rightarrow 0,
\end{equation*}
while $\langle\cdot,\cdot,\cdot\rangle:C\otimes C\otimes
C\rightarrow C[-1]$ is a chain homotopy
\begin{equation*}
    \langle\cdot,\cdot,\cdot\rangle:[\cdot,[\cdot,\cdot]]-[[\cdot,\cdot],\cdot]-[\cdot,[\cdot,\cdot]]*\sigma_{12}\Rightarrow 0
\end{equation*}
The coherence conditions satisfied by $S$ and $J$ translate to
equations involving $\langle\cdot,\cdot\rangle$ and
$\langle\cdot,\cdot,\cdot\rangle$, defining a structure which we
shall now describe.

\begin{defn}
A \emph{2-term $EL_\infty$-algebra} is a 2-term cochain complex $C$,
equipped with the following structure:
\begin{itemize}
\item a chain map $[\cdot,\cdot]:C\otimes C\rightarrow C$,
\item a chain homotopy $\langle\cdot,\cdot\rangle:[\cdot,\cdot]+[\cdot,\cdot]*\sigma\Rightarrow 0$,
\item a chain homotopy $\langle\cdot,\cdot,\cdot\rangle:[\cdot,[\cdot,\cdot]]-[[\cdot,\cdot],\cdot]-[\cdot,[\cdot,\cdot]]*\sigma_{12}\Rightarrow 0$,
\end{itemize}
such that the following equations hold:
\begin{equation}\label{eqn:l2l3=l3l2}
\begin{array}{c}
  [x,\langle y,z,w\rangle]+\langle x,[y,z],w\rangle+\langle x,z,[y,w]\rangle+[\langle x,y,z\rangle,w]+[z,\langle x,y,w\rangle]= \\
  =\langle x,y,[z,w]\rangle+\langle[x,y],z,w\rangle+[y,\langle x,z,w\rangle]+\langle y,[x,z],w\rangle+\langle y,z,[x,w]\rangle
\end{array}
\end{equation}
\begin{equation}\label{eqn:s12}
    \langle x,y,z\rangle+\langle y,x,z\rangle=-[\langle x,y\rangle,z]
\end{equation}
\begin{equation}\label{eqn:s23}
    \langle x,y,z\rangle+\langle x,z,y\rangle=[x,\langle y,z\rangle]-\langle[x,y],z\rangle-\langle y,[x,z]\rangle
\end{equation}
\begin{equation}\label{eqn:symm3}
    {\langle x,[y,z]\rangle} = \langle [y,z],x\rangle
\end{equation}
\end{defn}

\begin{rem}
We do not assume that $\langle\cdot,\cdot\rangle$ is symmetric, nor
that $\langle\cdot,\cdot,\cdot\rangle$ is skew-symmetric: in fact,
equations \eqref{eqn:s12}, \eqref{eqn:s23} and \eqref{eqn:symm3}
describe the symmetry properties of
$\langle\cdot,\cdot,\cdot\rangle$. These equations, however, are
easily seen to imply the following:
\begin{eqnarray}
  {[}\langle x,y\rangle,z] &=& [\langle y,x\rangle, z] \label{eqn:symm1}\\
  {[}x,\langle y,z\rangle] &=& [x,\langle z,y\rangle] \label{eqn:symm2}
\end{eqnarray}
In addition, equations \eqref{eqn:sym1}, \eqref{eqn:sym2} and
\eqref{eqn:sym3} obviously imply
\begin{eqnarray}
  d\langle x,y\rangle &=& d\langle y,x\rangle \label{eqn:symm4}\\
  \langle da,x\rangle &=& \langle x,da\rangle \label{eqn:symm5}
\end{eqnarray}
\end{rem}

Conversely, given a 2-term $EL_\infty$-algebra
$(C,[\cdot,\cdot],\langle\cdot,\cdot\rangle,\langle\cdot,\cdot,\cdot\rangle)$,
we can define a bracket on $L=\Gamma(C)$ by the formula
\eqref{eqn:catbracket}, with $[a,b]$ given by formula
\eqref{eqn:derivedbracket}, and $S$ and $J$ defined by formulas
\eqref{eqn:S} and \eqref{eqn:J}, respectively. The derived bracket
$[a,b]$ satisfies equations

\begin{eqnarray}
     [a,b]+[b,a]                  &=&\langle da,db\rangle \\
  {[}a,[b,c]]-[[a,b],c]-[b,[a,c]] &=& \langle da,db,dc\rangle
\end{eqnarray}
and it is routine to check that the axioms of a Lie 2-algebra are
satisfied.

Next, we apply normalization to a morphism $(F,F^2):L\rightarrow
L'$. The functor $F:L\rightarrow L'$ induces a chain map
$f=(f^0,f^1):C\rightarrow C'$, while $F^2$ can be written in the
form
$$F^2_{x,y}=([f^0(x),f^0(y)],-f^2(x,y)),$$
where $f^2:C^0\otimes C^0\rightarrow C^{-1}$ satisfies the following:
\begin{eqnarray}
  {[}f^0(x),f^0(y)]'-f^0([x,y]) &=& d'f^2(x,y) \\
  {[}f^1(a),f^0(y)]'-f^1([a,y]) &=& f^2(da,y) \\
  {[}f^0(x),f^1(b)]'-f^1([x,b]) &=& f^2(x,db)
\end{eqnarray}
In other words, $f^2$ is a homotopy from $[\cdot,\cdot]'*(f\otimes f)$ to $f*[\cdot,\cdot]$
\begin{defn}
A \emph{morphism} of $EL_\infty$-algebras from
$(C,[\cdot,\cdot],\langle\cdot,\cdot\rangle,\langle\cdot,\cdot,\cdot\rangle)$
to
$(C',[\cdot,\cdot]',\langle\cdot,\cdot\rangle',\langle\cdot,\cdot,\cdot\rangle')$
consists of
\begin{itemize}
\item a chain map $f=(f^0,f^1):C\rightarrow C'$
\item a chain homotopy $f^2:[\cdot,\cdot]'*(f\otimes f)\Rightarrow f*[\cdot,\cdot]$
\end{itemize}
such that the following equations hold:
\begin{equation}\label{eqn:f21}
    \langle f^0(x),f^0(y)\rangle'-f^1(\langle x,y\rangle)=f^2(x,y)+f^2(y,x)
\end{equation}
\begin{equation}\label{eqn:f22}
     \begin{array}{l}
       \langle f^0(x),f^0(y),f^0(z)\rangle'-f^1(\langle x,y,z\rangle)=\\
       =[f^0(x),f^2(y,z)]'-[f^0(y),f^2(x,z)]'-[f^2(x,y),f^0(z)]'-\\
       -f^2([x,y],z)-f^2(y,[x,z])+f^2(x,[y,z])
     \end{array}
\end{equation}
\end{defn}

\begin{defn}
The \emph{composite} of two morphisms $(f,f^2):C\rightarrow C'$ and
$(g,g^2):C'\rightarrow C''$ is defined to be $(gf,(g*f)^2)$, where
$$(g*f)^2(x,y)=g^2(f^0(x),f^0(y)))+g^1(f^2(x,y))$$
\end{defn}

Similarly, writing a 2-morphism $\Theta:(F,F^2)\Rightarrow(G,G^2)$
as
$$\Theta_x=(F(x),-\theta(x)),$$
with $\theta$ a chain homotopy from $f$ to $g$, leads to

\begin{defn}
A 2\emph{-morphism} $\theta:f\Rightarrow g$ is a chain homotopy
satisfying
$$f^2(x,y)-g^2(x,y)=[f^0(x),\theta(y)]'+[\theta(x),f^0(y)]'-\theta([x,y])-[\theta(x),\theta(y)]',$$
where $[\theta(x),\theta(y)]'$ is given by the formula
\eqref{eqn:derivedbracket}
\end{defn}

Conversely, morphisms and 2-morphisms of 2-term $EL_\infty$-algebras
induce the same of the corresponding Lie 2-algebras. We summarize
the above discussion in the following

\begin{thm}\label{thm:2catequiv}
2-term $EL_\infty$-algebras form a 2-category
$\mathbf{2TermEL_\infty}$ with the structure just defined. The
Dold-Kan correspondence induces an equivalence of 2-categories
\begin{diagram}
\mathbf{2Lie}&\pile{\rTo^N\\ \lTo_\Gamma}&\mathbf{2TermEL_\infty}
\end{diagram}
\end{thm}

\subsection{Special cases}

Setting the alternator $S$ to be the identity yields the notion of a
semistrict Lie 2-algebra. It coincides with the one defined in
\cite{BaezCrans}, since the bracket and the Jacobiator are (forced
to be) completely skew-symmetric in this case. Semi-strict Lie
2-algebras form a full sub-2-category corresponding to ordinary
2-term $L_\infty$-algebras upon normalization, as already shown in
\cite{BaezCrans}.

On the other hand, setting the Jacobiator $J$ to be the identity, we
get a full sub-2-category of hemistrict Lie 2-algebras. The
normalized complex $d:C^{-1}\longrightarrow C^0$ of such a 2-algebra
inherits the structure of a differential graded \emph{Leibniz}
algebra, since the right hand sides of the equations
\ref{eqn:Leib1}--\ref{eqn:Leib4} vanish. In particular, $C^0$ is a
Leibniz algebra acting on $C^{-1}$ on both sides. $C^{-1}$ itself
becomes a Leibniz algebra with respect to the derived bracket
(defined, as usual, by formula \eqref{eqn:derivedbracket}), making
$d:C^{-1}\longrightarrow C^0$ a Leibniz algebra crossed module.
Representations of Leibniz algebras and crossed modules were
considered in \cite{LodPir}

But this is not all. In addition, we have a bilinear operation
\begin{equation*}
    \langle\cdot,\cdot\rangle:C^0\otimes C^0\rightarrow C^{-1}
\end{equation*}
which measures, via equations \eqref{eqn:sym1}, \eqref{eqn:sym2} and
\eqref{eqn:sym3}, the failure of the Leibniz algebra $C^0$ to be a
Lie algebra, as well as the failure of the representation $C^{-1}$
of $C^0$ to be symmetric (in the terminology of \cite{LodPir}). It
obeys the equations

\begin{eqnarray*}
  [\langle x,y\rangle,z] &=& 0 \\
  {[}x,\langle y,z\rangle] &=& \langle[x,y],z\rangle+\langle y,[x,z]\rangle  \\
  \langle x,[y,z]\rangle &=& \langle[y,z],x\rangle
\end{eqnarray*}

Equations \eqref{eqn:symm2}, \eqref{eqn:symm4} and \eqref{eqn:symm5}
are implied, as in the general case. In particular, the image of
$\langle\cdot,\cdot\rangle$ is an anti-symmetric submodule of
$C^{-1}$ (in the terminology of \cite{LodPir}), while the
skew-symmetric part of $\langle\cdot,\cdot\rangle$ is annihilated by
the action of $C^0$ on both sides and is contained in the kernel of
$d$.

\begin{example}
Given a Leibniz algebra $\mathfrak{g}$, denote by
$\mathfrak{g}^{\mathrm{ann}}$ the subspace of $\mathfrak{g}$ spanned
by the elements of the form $[x,x]$, $x\in\mathfrak{g}$. It is in
fact a two-sided ideal in $\mathfrak{g}$. Setting
$C^0=\mathfrak{g}$, $C^{-1}=\mathfrak{g}^{\mathrm{ann}}$ and $d$ the
inclusion map, we get a dg Leibniz algebra. Moreover, setting
\begin{equation*}
    \langle x,y\rangle=[x,y]+[y,x]
\end{equation*}
gives an alternator with all the required properties. Thus,
\emph{any} Leibniz algebra gives rise to a hemistrict Lie 2-algebra,
albeit a rather special one: the alternator is symmetric, and
$C^{-1}$ is an anti-symmetric $C^0$-module. However, if the
characteristic of the ground field is different from 2, any
hemistrict Lie 2-algebra with surjective $\langle\cdot,\cdot\rangle$
and injective $d$ is of this form.
\end{example}

\begin{example}\label{eg:quadratichs} Let $\mathfrak{g}$ be a Lie algebra equipped with an
$ad$-invariant symmetric bilinear form $\langle\cdot,\cdot\rangle$.
Setting $C^0=\mathfrak{g}$, $C^{-1}$=$\mathbf{k}$, $d=0$ gives rise
to a hemi-strict Lie 2-algebra with $[x,a]=-[a,x]=0$ and the
alternator given by $\langle\cdot,\cdot\rangle$.
\end{example}

Finally, if both the alternator and the Jacobiator are trivial (i.e.
the Lie 2-algebra is strict), we get a differential graded Lie
algebra on the normalized complex. The derived bracket
\eqref{eqn:derivedbracket} is then a Lie bracket, yielding a Lie
algebra crossed module.

\section{Skew-symmetrization}\label{sec:Skew-sym}

In this section we assume that the characteristic of the ground
field $\mathbf{k}$ is different from 2 or 3.

Suppose $(V,[\cdot,\cdot],S,J)$ is a Lie 2-algebra,
$(C,[\cdot,\cdot],\langle\cdot,\cdot\rangle,\langle,\cdot,\cdot,\cdot\rangle)$
-- the corresponding $EL_\infty$-algebra. Define multilinear
skew-symmetric maps $\{\cdot,\cdot\}:\bigwedge^2C\rightarrow C$ and
$\{\cdot,\cdot,\cdot\}:\bigwedge^3C\rightarrow C[-1]$ as follows:
\begin{eqnarray*}
  \{x,y\} &=& \frac{1}{2}([x,y]-[y,x]) \\
  \{x,a\} &=& \frac{1}{2}([x,a]-[a,x])=-\{a,x\} \\
  \{x,y,z\} &=& [x,y,z]-T(x,y,z),
\end{eqnarray*}
where
\begin{eqnarray*}
  [x_1,x_2,x_3] &=& \frac{1}{6}\sum_{\sigma\in S_3}(-1)^\sigma\langle x_{\sigma(1)},x_{\sigma(2)},x_{\sigma(3)}\rangle \\
  T(x_1,x_2,x_3) &=& \frac{1}{12}\sum_{\sigma\in S_3}(-1)^\sigma\langle[x_{\sigma(1)},x_{\sigma(2)}],x_{\sigma(3)}\rangle
\end{eqnarray*}
Then we can prove the following
\begin{prop}\label{prop:skewsym}
$(C,\{\cdot,\cdot\},\{\cdot,\cdot,\cdot\})$ is an
$L_\infty$-algebra.
\end{prop}
which can be expanded to
\begin{thm}\label{thm:skewsym}
Skew-symmetrization defines a projection 2-functor
$$SS:\mathbf{2Lie}\rightarrow \mathbf{SS2Lie}$$ onto the 2-category of
semistrict Lie 2-algebras.
\end{thm}
The proof is a routine verification of the axioms.
\begin{cor}
In particular, there is a skew-symmetrization 2-functor
$$SS:\mathbf{HS2Lie}\rightarrow \mathbf{SS2Lie}$$
from hemistrict to semistrict Lie 2-algebras with
\begin{equation*}
    \{x,y,z\} = -T(x,y,z)
\end{equation*}
\end{cor}

\begin{example}\label{eg:quadraticss}
Applying this to the hemistrict Lie 2-algebra of Example
\ref{eg:quadratichs} yields a semistrict Lie 2-algebra with the same
underlying category and brackets, with the Jacobiator given by
$$\langle\cdot,\cdot,\cdot\rangle=-\frac{1}{2}\langle[\cdot,\cdot],\cdot\rangle$$
In the case when $\mathfrak{g}$ is semisimple, with
$k=\langle\cdot,\cdot\rangle$ the Killing form, the resulting Lie
2-algebra is the string Lie 2-algebra denoted by
$\mathfrak{g}_{-\frac{1}{2}}$ in \cite{BaezCrans} and
$\mathfrak{str}_{-\frac{1}{2}}(\mathfrak{g})$ in \cite{Hen}.
\end{example}

% Apply this to a Leibniz algebra

\section{Categorical and Homotopy invariance}\label{sec:Homotopy}

Recall that an equivalence of linear categories consists of a pair
of linear functors $F:V\rightarrow V'$, $G:V'\rightarrow V$,
together with linear natural transformations $\Phi:F*G\Rightarrow
1_V'$ and $\Psi:G*F\Rightarrow 1_V$. It is a standard fact in
category theory that a functor $F$ induces an equivalence of
categories if and only if it is fully faithful and essentially
surjective; its quasi-inverse $G$ is then unique up to natural
isomorphism. This carries over to the linear case with obvious
modifications.

Using the Dold-Kan correspondence, it is easy to deduce that a
linear functor $F$ is fully faithful and essentially surjective if
and only if $f=N(F)$ is a quasi-isomorphism (a chain map inducing
isomorphism in cohomology), and if $G$ is a quasi-inverse to $F$,
then $g=N(G)$ is a homotopy inverse to $f$, and vice versa. In
particular, if $C$ is a cochain complex of vector spaces and $H$ is
its cohomology, viewed as a complex with zero differential, then
there exists a homotopy equivalence $C\pile{\rightarrow\\
\leftarrow}H$ (a Hodge decomposition).

A categorically invariant algebraic structure is, heuristically, a
structure that can be transferred along categorical equivalences.
This means that if $V$ is a category equipped with this structure
and $F:V\rightarrow V'$ is fully faithful and essentially
surjective, then there exists the same type of structure on $V'$,
unique up to equivalence, such that $F$ induces an equivalence of
categories with the structure. Making this precise in full
generality requires a categorification of the notion of operad. We
shall not attempt this here, taking advantage instead of the
Dold-Kan correspondence in order to transfer everything to chain
complexes, where the similar notion of a homotopy-invariant
algebraic structure is well-known.

Recall also that a category is called \emph{skeletal} if isomorphic
objects are equal. It follows that a linear category is skeletal if
and only if its normalized complex has zero differential. By Hodge
decomposition, every linear category is equivalent to a skeletal
one. We call a Lie 2-algebra skeletal if its underlying linear
category is.

We have the following result.

\begin{thm}\label{thm:invariance}
The structure of 2-term $EL_\infty$-algebra is homotopy-invariant;
equivalently, the structure of a Lie 2-algebra is categorically
invariant. In particular, every Lie 2-algebra is equivalent, as a
Lie 2-algebra, to a skeletal one.
\end{thm}
The proof is a standard exercise in homological perturbation theory.
\begin{rem}
That semistrict Lie 2-algebras are categorically invariant follows
from the well-known fact that $L_\infty$-algebras are
homotopy-invariant. However, hemistrict Lie 2-algebras are
\emph{not} categorically invariant.
\end{rem}
It remains to determine what skeletal Lie 2-algebras look like,
and to classify them up to equivalence.

So, let $L$ be a skeletal Lie 2-algebra, $C$ its normalized complex,
with $d=0$. Because of this last fact, $C^0$ is an honest Lie
algebra, acting on $C^{-1}$, with $[x,a]=-[a,x]$. The Jacobiator
$\langle\cdot,\cdot,\cdot\rangle$, which we shall here rename $j$,
obeys the equation \eqref{eqn:l2l3=l3l2}, which can be rewritten in
the form
\begin{equation}\label{eqn:3cocycle1}
    \begin{array}{l}
       [x,j(y,z,w)]-[y,j(x,z,w)]+[z,j(x,y,w)]+[j(x,y,z),w]- \\
       -j([x,y],z,w)-j(y,[x,z],w)-j(y,z,[x,w])+ \\
       +j(x,[y,z],w)+j(x,z,[y,w])-j(x,y,[z,w])=0
     \end{array}
\end{equation}
The reader can recognize this equation as saying that $j$ is a
3-cocycle in the Loday-Pirashvili complex for $C^0$, viewed as a
Leibniz algebra, with coefficients in the (symmetric) representation
$C^{-1}$, as defined in \cite{LodPir}. It is not, however, a
Chevalley-Eilenberg cocycle, for lack of skew-symemtry. In fact,
with the alternator denoted by $s$, equations \eqref{eqn:s12},
\eqref{eqn:s23} and \eqref{eqn:symm3} translate in this case to
\begin{eqnarray}
  j(x,y,z)+j(y,x,z) &=& [z,s(x,y)]\label{eqn:3cocycle2} \\
  j(x,y,z)+j(x,z,y) &=& [x,s(y,z)]-s([x,y],z)-s(y,[x,z])\label{eqn:3cocycle3} \\
  s([x,y],z) &=& s(z,[x,y])\label{eqn:3cocycle4}
\end{eqnarray}
To see when two skeletal Lie 2-algebras are equivalent, we first
remark that skeletal categories are equivalent if and only if they
are strictly isomorphic, and that a morphism between skeletal Lie
2-algebras is, in particular, a strict homomorphism of the
underlying Lie algebras and representations. Because of this, it is
necessary and sufficient to determine when
$$(1,f):(C,[\cdot,\cdot],s,j)\longrightarrow(C,[\cdot,\cdot],s',j')$$
is a morphism of $EL_\infty$-algebras, with $f:C^0\otimes
C^0\rightarrow C^{-1}$. The condition is given by equations
\eqref{eqn:f21} and \eqref{eqn:f22} which read in this case:
\begin{eqnarray}
  s'(x,y)-s(x,y) &=& f(x,y)+f(y,x) \label{eqn:3cob1}\\
  j'(x,y,z)-j(x,y,z) &=& \begin{array}{c}
                           [x,f(y,z)]-[y,f(x,z)]-[f(x,y),z]- \\
                           -f([x,y],z)-f(y,[x,z])+f(x,[y,z])
                         \end{array}\label{eqn:3cob2}
\end{eqnarray}
The second equation says that Loday-Pirashvili cocycles $j'$ and $j$
are cohomologous. Setting
$$g(x,y)=-f(x,y)+[x,\theta(y)]+[\theta(x),y]-\theta([x,y])$$
with an arbitrary $\theta:C^0\rightarrow C^{-1}$ then gives a
morphism $(1,g)$ the other way, with $\theta:(1,g)*(1,f)\Rightarrow
(1,0)$.

In general, given a Lie algebra $\mathfrak{g}$ and a representation
$M$, define $ZL^3_{\mathrm{Lie}}(\mathfrak{g},M)$ to be the space
consisting of pairs $(s,j)$, where
$s:\mathfrak{g}^{\otimes^2}\rightarrow M$ and
$j:\mathfrak{g}^{\otimes^3}\rightarrow M$ satisfy equations
(\ref{eqn:3cocycle1}-\ref{eqn:3cocycle4}); define
$BL^3_{\mathrm{Lie}}(\mathfrak{g},M)$ to be the space of pairs
$(s,j)$ such that there exists an
$f:\mathfrak{g}^{\otimes^2}\rightarrow M$ making equations
\eqref{eqn:3cob1} and \eqref{eqn:3cob2} hold with $s$ (resp. $j$) on
the left hand side. We can prove an easy

\begin{lemma}
$BL^3_{\mathrm{Lie}}(\mathfrak{g},M)\subseteq
ZL^3_{\mathrm{Lie}}(\mathfrak{g},M)$
\end{lemma}

\begin{rem}
The proof of the lemma depends not only on the equations
(\ref{eqn:Leib1}-\ref{eqn:Leib4}) (with vanishing right hand sides),
as for the differential in the Loday-Pirashvili complex, but also on
the skew-symmetry
\begin{eqnarray*}
  [x,y] &=& -[y,x] \\
  {[}x,a] &=& -[a,x]
\end{eqnarray*}
Lie algebras in this context are to be viewed as Leibniz algebras
which are also skew-symmetric.
\end{rem}
We define
$HL^3_{\mathrm{Lie}}(\mathfrak{g},M)=ZL^3_{\mathrm{Lie}}(\mathfrak{g},M)/BL^3_{\mathrm{Lie}}(\mathfrak{g},M)$.
With this we have obtained the following
\begin{thm}\label{thm:classification}
Skeletal Lie 2-algebras are classified up to equivalence by the
following data:
\begin{itemize}
\item a Lie algebra $\mathfrak{g}$
\item a representation $M$ of $\mathfrak{g}$
\item a class $[(s,j)]\in HL^3_{\mathrm{Lie}}(\mathfrak{g},M)$
\end{itemize}
\end{thm}

\begin{rem}
This classification generalized the classification of semistrict
skeletal Lie 2-algebras due to Baez and Crans \cite{BaezCrans}, for
our cohomology space then reduces to the Chevalley-Eilenberg
cohomology. But it is also remarkably similar to the classification
of skeletal braided categorical groups, due to Joyal and Street
(Proposition 3.1 of \cite{JS}). In that classification, a very
similar construction, due to Eilenberg and MacLane and going back to
1950, was used to obtain a group $H^3_{\mathrm{ab}}(G,M)$ for an
abelian group $G$ with coefficients in a representation $M$. We
suspect that Lie 2-algebras may be related to linear braided
monoidal categories by a a sort of categorified Koszul duality, yet
to be described.
\end{rem}

\begin{thm}\label{thm:LPtoChE}
The skew-symmetrization functor induces a map
$$\mathrm{ss}:HL^3_{\mathrm{Lie}}(\mathfrak{g},M)\longrightarrow
H^3(\mathfrak{g},M)$$ onto the Chevalley-Eilenberg cohomology,
fitting into the exact sequence
\begin{equation*}
    0\longrightarrow \mathrm{Hom}(\wedge^2\mathfrak{a},M)\stackrel{\iota}{\longrightarrow}
    HL^3_{\mathrm{Lie}}(\mathfrak{g},M)\stackrel{\mathrm{ss}}{\longrightarrow}H^3(\mathfrak{g},M)\longrightarrow
    0
\end{equation*}
where
$$\mathfrak{a}=\mathfrak{g}/[\mathfrak{g},\mathfrak{g}]$$
is the abelianization, and $\iota(a)=[(a,0)]$. A canonical splitting
is given by $\phi\mapsto[(0,\phi)]$.
\end{thm}
\begin{example}
Let $\mathfrak{g}$ be a semisimple Lie algebra, with Killing form
$k$ and Cartan tensor $\phi=k([\cdot,\cdot],\cdot)$. Since
semi-simple Lie algebras are perfect, $\mathfrak{a}=0$, hence
$\mathrm{ss}$ is an isomorphism. It sends the class of $(k,0)$ to
that of $-\frac{1}{2}\phi$ (Example \ref{eg:quadraticss}). Notice
that $\phi$ is the Loday-Pirashvili coboundary of $k$, hence $(k,0)$
is cohomologous to $(0,-\frac{1}{2}\phi)$. The corresponding
hemistrict and semistrict Lie 2-algebras are equivalent as Lie
2-algebras and, since the cohomology space is one-dimensional, any
Lie 2-algebra structure with this underlying Lie algebra and module
is equivalent to a multiple of either.
\end{example}
% Add example of quadratic Lie algebras here

\section{Applications}\label{sec:Appl}

In this section we assume that the characteristic of the ground
field $\mathbf{k}$ is zero.

Recall that, in deformation theory, one considers
$L_\infty$-algebras of the form

$$(L=\bigoplus_{k\in(-n,+\infty)}L^k,\delta=[\cdot],[\cdot,\cdot],[\cdot,\cdot,\cdot],\ldots)$$
where the $k$-nary bracket is of degree $2-k$ and $n$ is a
non-negative integer. One is then interested in the space of
solutions of the (generalized) Maurer-Cartan equation

\begin{equation*}
    \mathrm{MC}(L)=\{\gamma\in L^1|\sum_{k\geq1}{\frac{1}{k!}}[\gamma,\ldots,\gamma]=0\}
\end{equation*}
where the $k$th summand is the $k$-nary bracket of $\gamma$ with
itself (in particular, $[\gamma]=\delta(\gamma)$ is the
differential). The equation describes the type of algebraic
structure one wishes to study. To make sense of it, one either
renders $L$ nilpotent by tensoring with a nilpotent commutative
algebra, or assumes that all but a finite number of $k$-nary
brackets are zero; in fact, in most cases occurring in practice only
the differential and the binary bracket are nontrivial, making $L$ a
differential graded Lie algebra (dgla).

There is an equivalence relation on $\mathrm{MC}(L)$ induced by the
infinitesimal action of $L^0\ni x\mapsto\bar{x}$ where
$$\bar{x}(\gamma)=\delta(x)+[\gamma,x]+\frac{1}{2}[\gamma,\gamma,x]+\cdots$$
If $L$ is a dgla with no components of negative degree ($n=1$),
$L^0$ is a Lie algebra; integrating its action (for nilpotent $L$)
gives rise to an action groupoid, known as the Deligne groupoid of
$L$. This groupoid presents the moduli stack of $L$, which is the
main object of study in deformation theory.

However, the presence of negative degrees leads to a richer
structure involving higher symmetries. This was first noticed by
Deligne who constructed, in an unpublished letter to L. Breen, a
strict 2-groupoid over $\mathrm{MC}(L)$ for a dgla $L$ with $n=2$.
This construction was rediscovered by Getzler \cite{Get02}. It uses
the derived bracket on $L^{-1}$ parametrized by
$\gamma\in\mathrm{MC}(L)$. In particular, the 2-group of
automorphisms of $\gamma$ is obtained by integrating the
corresponding Lie algebra crossed module under the derived bracket
(i.e. a strict Lie 2-algebra in our sense).

It was Getzler \cite{Get04} who generalized this construction to
$L_\infty$-algebras and higher values of $n$. The result is a weak
$n$-groupoid which, by definition, is a Kan complex with unique
fillers for horns in dimension higher than $n$. However, much
remains to be understood about the structure of higher symmetries
even at the infinitesimal level. Here we present a construction, for
a dgla with $n=3$, of a kind of crossed module involving Lie
2-algebras of our kind.

To begin, notice that there is a family of $L_\infty$-algebras
parametrized by $\mathrm{MC}(L)$:
\begin{eqnarray*}
  \delta_\gamma &=& \delta+[\gamma,\cdot]+\frac{1}{2}[\gamma,\gamma,\cdot]+\cdots \\
  {[\cdot,\cdot]_\gamma} &=& [\cdot,\cdot]+[\gamma,\cdot,\cdot]+\frac{1}{2}[\gamma,\gamma,\cdot,\cdot]+\cdots \\
  {[\cdot,\cdot,\cdot]_\gamma} &=& [\cdot,\cdot,\cdot]+[\gamma,\cdot,\cdot,\cdot]+\cdots \\
  \cdots && \cdots
\end{eqnarray*}
and that truncation
$$L^{-n+1}\stackrel{\delta_\gamma}{\rightarrow}\cdots\stackrel{\delta_\gamma}{\rightarrow}L^{-1}\stackrel{\delta_\gamma}{\rightarrow}\bar{L}^0$$
where $\bar{L}^0=\ker\delta_\gamma$ defines an
$L_\infty$-subalgebra, the algebra of infinitesimal automorphisms of
$\gamma$. It is this algebra that we shall presently study.

\subsection{Case $n=2$} Here we have a 2-term $L_\infty$-algebra of
the form
$$L^{-1}\stackrel{d}{\longrightarrow}\bar{L}^0$$
with only $d=\delta_\gamma$, $[\cdot,\cdot]_\gamma$ and
$[\cdot,\cdot,\cdot]_\gamma$ nontrivial. This gives rise to a
semistrict Lie 2-algebra, with the derived bracket, defined by
formula \eqref{eqn:derivedbracket}, giving the crossed module
structure.

\begin{example}
Let $M$ be a smooth manifold. Set $L=\Gamma(\wedge^\cdot TM)[1]$,
the algebra of smooth multivector fields. It is a dgla under the
Schouten bracket and zero differential. A solution $\gamma$ of the
Maurer-Cartan equation is, by definition, a Poisson structure on
$M$. The strict Lie 2-algebra of infinitesimal automorphisms of
$\gamma$ has $L^{-1}=C^\infty(M)$, $\bar{L}^0$ the space of Poisson
vector fields, and $d=\delta_\gamma=[\gamma,\cdot]$ the Lichnerowicz
differential, sending a function $a$ to its Hamiltonian vector
field. The derived bracket on $L^{-1}$ is just the Poisson bracket
of functions determined by $\gamma$. This is the main example whose
integration was given in \cite{Get02}. It can be generalized to any
Lie bialgebroid.
\end{example}

\begin{example}
With $L$ as above, let $H$ be a 3-form on $M$. It extends by the
Leibniz rule to define a trilinear operation $[\cdot,\cdot,\cdot]$
on $L$ of degree $-1$. Together with the Schouten bracket, it
defines an $L_\infty$-structure if and only if $H$ is closed. A
solution $\gamma$ of the Maurer-Cartan equation is an $H$-twisted
Poisson structure (\cite{SevWe}, \cite{Roy3-QuasiLie}). The
construction of the Lie 2-algebra of infinitesimal automorphism
proceeds as in the previous example, except now it is only
semistrict. This example generalizes to any quasi-Lie bialgebroid
\cite{Roy3-QuasiLie}.
\end{example}

\subsection{Case $n=3$, dgla.} In this case the truncated dgla is 3-term:

$$L^{-2}\stackrel{d}{\longrightarrow}L^{-1}\stackrel{d}{\longrightarrow}\bar{L}^0$$
with $d=\delta_\gamma=\delta+\{\gamma,\cdot\}$ and $\{\cdot,\cdot\}$
the bracket on $L$. In particular, $\bar{L}^0$ is a Lie algebra
acting on $L^{-1}$ and $L^{-2}$ in a way compatible with $d$, but
there is also a \emph{symmetric} bilinear map

$$\{\cdot,\cdot\}:L^{-1}\otimes L^{-1}\longrightarrow L^{-2}$$
which will play the role of an alternator, so let us denote it by
$\langle\cdot,\cdot\rangle$ from now on. Set $C^{i}=L^{i-1}$,
$i=-1,0$, and introduce the derived brackets on $C$ as follows:
\begin{eqnarray*}
  {[x,y]} &=& \{dx,y\} \\
  {[x,a]} &=& \{dx,a\} \\
  {[a,x]} &=& 0
\end{eqnarray*}
for $x,y\in C^0$, $a\in C^{-1}$.

Furthermore, viewing $\bar{L}^0$ as a cochain complex concentrated
in degree 0 or, equivalently, as a linear category with only
identity arrows, we get a chain map
$$(d,0):(C^0,C^{-1})\longrightarrow (\bar{L}^0,\{0\}),$$
so that $\partial=\Gamma(d,0):\Gamma(C)\longrightarrow\bar{L}^0$ is
a linear functor.

Lastly, since the action of $\bar{L}^0$ (given by
$T\mapsto\{T,\cdot\}=-\{\cdot,T\}$) commutes with $d$, it induces a
functorial action of $\bar{L}^0$ on the linear category $\Gamma(C)$.
We have the following

\begin{thm}\label{thm:2crossed}
For any dgla $L$ concentrated in degrees $(-3,0]$,
\begin{itemize}
\item $\Gamma(C)$ is a hemistrict Lie 2-algebra;
\item $\partial:\Gamma(C)\longrightarrow\bar{L}^0$ is a morphism of Lie
2-algebras;
\item $\bar{L}^0$ acts on $\Gamma(C)$ by strict derivations.
\end{itemize}
In addition, the following crossed module identities hold:
\begin{eqnarray*}
  \partial\{T,\mathbf{f}\} &=& \{T,\partial\mathbf{f}\} \\
  \{\partial\mathbf{f},\mathbf{g}\} &=& [\mathbf{f},\mathbf{g}]
\end{eqnarray*}
where $T\in\bar{L}^0$, $\mathbf{f},\mathbf{g}\in \Gamma(C)$.
\end{thm}
The proof is a routine verification.
%expand this example to a subsection.
\begin{example}
Let $E\rightarrow M$ be a vector bundle with a fiberwise smooth
inner product $\langle\cdot,\cdot\rangle$. In \cite{Roy4-GrSymp} we
constructed a dgla $L$ as above, such that the solutions of the
Maurer-Cartan equation are precisely Courant algebroid structures on
$E$, with the Courant bracket defined as the derived bracket.
Theorems \ref{thm:2crossed} and \ref{thm:skewsym} combine to yield
the main result of \cite{RoyWe} as an immediate corollary.
\end{example}

%\section{Concluding remarks and further questions}\label{sec:Concl}

%\bibliographystyle{plain}
%\bibliography{../ref}

\end{document}